# Power spectrum of the fluctuation of Chebyshev's prime counting function

Boon Leong Lan and Shaohen Yong

School of Engineering, Monash University, 46150 PJ, Selangor, Malaysia

**Abstract**

The one-sided power spectrum of the fluctuation of Chebyshev's weighted prime counting function is numerically estimated based on samples of the fluctuating function of different sizes. The power spectrum is also estimated analytically for large frequency based on Riemann hypothesis and the exact formula for the fluctuating function in terms of all the non-trivial Riemann zeroes. Our analytical estimate is consistent with our numerical estimate of a $1/f^2$ power spectrum.

*PACS*: 05.45.Tp

*Keywords*: Time series; Fluctuation; Chebyshev's prime counting function; Power-law power spectrum; Riemann hypothesis

The studies of prime numbers using the tools of statistical physics have proven to be quite fruitful. Wolf [1] found that the primes have multifractal properties. He [2] also found that the 'time' series formed using the number of primes in successive intervals of equal length $l$ has the same power-law ($1/f^{1.64}$) power spectrum for different $l$. Regarding $x$ as a 'time' variable, Gamba et al. [3] found that the fluctuation

$$R(x) - \pi(x) \tag{1}$$

of the prime counting function $\pi(x)$ from its smooth average part given by Riemann's function $R(x)$ has a positive Liapunov exponent. The prime counting function

$$\pi(x) = \sum_{p \leq x} 1 \tag{2}$$

is a staircase function with a jump of one when $x$ crosses a prime $p$; at $x=p$, the function takes a value halfway between its new and old values. Our [4] recent numerical estimates suggest that the fluctuation $R(x) - \pi(x)$ has a $1/f^2$ power spectrum.

In this paper, we study, both numerically and analytically, the power spectrum of the fluctuation

$$\psi_{fluc}(x) \equiv \psi(x) - \left( x + \sum_n \frac{x^{-2n}}{2n} - \ln 2\pi \right) \tag{3}$$

of Chebyshev's *weighted* prime counting function $\psi(x)$ from its smooth average part given by

$$x + \sum_n \frac{x^{-2n}}{2n} - \ln 2\pi .$$

Chebyshev's weighted prime counting function is [5] defined as

$$\psi(x) = \sum_{p^m \leq x} \ln p \tag{4}$$

where the sum is over prime powers. Chebyshev's weighted prime counting function $\psi(x)$ is [5] also a staircase function like the non-weighted prime counting function $\pi(x)$ but it has a jump of $\ln p$ when $x$ crosses a prime power $p^m$; at $x=p^m$, the function takes a value halfway between its new and old values. For e.g.,

$$\psi(10) = 3\ln 2 + 2\ln 3 + \ln 5 + \ln 7 = \ln 2520,$$

$$\psi(12) = \psi(10) + \ln(11) = \ln 27720,$$

$$\psi(11) = \frac{\psi(10) + \psi(12)}{2}.$$

The prime counting function $\psi(x)$, like $\pi(x)$, has [5] an exact formula

$$\psi(x) = x + \sum_n \frac{x^{-2n}}{2n} - \ln 2\pi - \sum_\rho \frac{x^\rho}{\rho} \tag{5}$$

in terms of all the nontrivial complex zeroes $\rho$ (they come in complex conjugate pairs) of the celebrated Riemann zeta function and is therefore of great importance in mathematics.

First we present our numerical result for the power spectrum of the fluctuating function $\psi_{fluc}(x)$. The one-sided power spectrum $P(f)$ of the fluctuating function $\psi_{fluc}(x)$ in Eq. (3) was numerically estimated based on different size samples of the function at positive integers $x$ starting at $x=2$ with a sampling interval of one. The fluctuating function $\psi_{fluc}(x)$ at each integer $x$ from $x=2$ to $x=10^4$ is shown Fig. 1. The power spectrum was, as usual [6], estimated as a function of frequency $f$ (frequency in arbitrary frequency unit since *x* is in arbitrary 'time' unit) up to the Nyquist frequency. Since the sampling interval $\Delta x = 1$, the Nyquist frequency, which is [6] defined as $1/(2\Delta x)$, is 0.5. The power spectrum of $\psi_{fluc}(x)$ was estimated using the maximum entropy (or all poles) method [6], and Welch's modified periodogram method, which uses fast fourier transform, in MATLAB. These two methods gave consistent power spectrum estimates in all cases and so only the maximum-entropy estimated power spectrum is presented. In general [6], in the maximum entropy method, if the number of poles used is too big, the estimated power spectrum will have spurious oscillations. For the fluctuating function $\psi_{fluc}(x)$, the power spectrum estimated using 2 poles already starts to exhibit such oscillations for frequency greater than 0.1. So we will present only the power spectrum estimated using 1 pole.

Fig. 2 shows the one-sided power spectrum of the fluctuating function $\psi_{fluc}(x)$ that was numerically estimated based on different size samples of the function: $10^3$, $10^4$, $10^5$, and $10^6$ sample points. As the size of the samples increases from $10^3$ to $10^5$, the estimated power spectrum converges to a straight line on the log-log scale. In all cases, the curve in the tail of the power spectrum near the Nyquist frequency is due to aliasing [6], i.e., power outside the frequency range from 0 to the Nyquist frequency is falsely folded over or aliased into the range. For the $10^5$ and $10^6$ cases, each of the estimated power spectrum is well-fitted by a straight line in the frequency range $10^{-3}$ to $10^{-1}$ and thus behaves as a power law

$$P(f) = af^b \qquad (6)$$

where *a* and *b* are constants. The fit yields *b* = -2 in each case. Fig. 3 shows the power spectrum of $\psi_{fluc}(x)$, which was estimated based on $10^6$ sample points, together with the excellent straight-line fit. Our numerical results therefore suggest that the one-sided power spectrum $P(f)$ of the fluctuation $\psi_{fluc}(x)$ of Chebyshev's weighted prime counting function is $1/f^2$.

Next we present our analytical result for the power spectrum of the fluctuating function $\psi_{fluc}(x)$. The one-sided power spectrum $P(f)$ of the fluctuating function $\psi_{fluc}(x)$ was estimated analytically, for large frequency $f$, based on the exact formula for $\psi_{fluc}(x)$ from Eq. (5) in terms of all the non-trivial Riemann zeroes with the assumption (Riemann hypothesis) that the real part of all the zeroes is one half:

$$\psi_{fluc}(x) = -\left( \sum_t \frac{x^{\frac{1}{2}} e^{it \ln x}}{\frac{1}{2} + it} + \sum_t \frac{x^{\frac{1}{2}} e^{-it \ln x}}{\frac{1}{2} - it} \right). \tag{7}$$

From Eq. (7), the Fourier transform $F$ of $\psi_{fluc}(x)$, $H(f)$, is given by

$$H(f) = -\left( \sum_t \frac{F\left\{ x^{\frac{1}{2}} e^{it \ln x} \right\}}{\frac{1}{2} + it} + \sum_t \frac{F\left\{ x^{\frac{1}{2}} e^{-it \ln x} \right\}}{\frac{1}{2} - it} \right). \tag{8}$$

First, consider the first term in Eq. (8). For large $f$, the Fourier transform of $x^{\frac{1}{2}} e^{it \ln x}$ can be evaluated approximately using the method [7] of stationary phase:

$$F\left\{ x^{\frac{1}{2}} e^{it \ln x} \right\} = \int dx\, x^{\frac{1}{2}}\, e^{i(t \ln x - fx)}$$

$$\approx \frac{\sqrt{-2\pi i}}{f^{3/2}} t\, e^{i\left( t \ln \frac{t}{f} - t \right)} \qquad f > 0. \tag{9}$$

The discrete sum over the Riemann zeroes is then replaced by an integral

$$\sum_t \to \int dt\, d(t)$$

where $d(t)$ is the density of Riemann zeroes at height $t$. To do the $t$ integral, $d(t)$ is approximated by its average [8] given by $\frac{1}{2\pi}\ln\left(\frac{t}{2\pi}\right)$, and for large $f$, the $t$ integral can also be evaluated approximately using the method of stationary phase yielding

$$\sum_t \frac{F\left\{x^{\frac{1}{2}}e^{it\ln x}\right\}}{\frac{1}{2}+it} \approx -\frac{\ln(f/2\pi)}{\frac{1}{2}+if}e^{-if} \qquad f>0. \tag{10}$$

Similar calculations for the second term in Eq. (8) yield

$$\sum_t \frac{F\left\{x^{\frac{1}{2}}e^{-it\ln x}\right\}}{\frac{1}{2}-it} \approx -\frac{\ln(|f|/2\pi)}{\frac{1}{2}-i|f|}e^{i|f|} \qquad f<0. \tag{11}$$

The one-sided power spectrum $P(f)$ of $\psi_{fluc}(x)$ is defined [6] as

$$P(f) \equiv |H(f)|^2 + |H(-f)|^2 \qquad f \geq 0. \tag{12}$$

Since $\psi_{fluc}(x)$ is real, $|H(f)|^2 = |H(-f)|^2$ [6], and so

$$P(f) = 2|H(f)|^2 \qquad f \geq 0. \tag{13}$$

For large $f$, Eq.(10) gives

$$P(f) \approx \frac{2\ln^2(f/2\pi)}{f^2}, \tag{14}$$

which is approximately $\frac{const.}{f^2}$.

Our analytical estimate of the one-sided power spectrum $P(f)$ of the fluctuating function $\psi_{fluc}(x)$ for large frequency $f$, which is based on Riemann hypothesis and the exact formula Eq. (7) for $\psi_{fluc}(x)$ in terms of all the non-trivial Riemann zeroes, is therefore consistent with our numerical estimate of $1/f^2$ based on Eq. (3) for $\psi_{fluc}(x)$ which involves only the prime numbers. This result adds to the evidence [9] in favor of Riemann hypothesis.

**Acknowledgement**

We sincerely thank the referee for very kindly pointing out how the analytical calculations of the Fourier transform and the summation in Eq. (8) could be done approximately.

**Fig. 1**

The fluctuating function $\psi_{fluc}(x)$, defined in Eq. (3), plotted at each integer $x$ from $x=2$ to $x=10^4$. The data points are connected by straight lines.

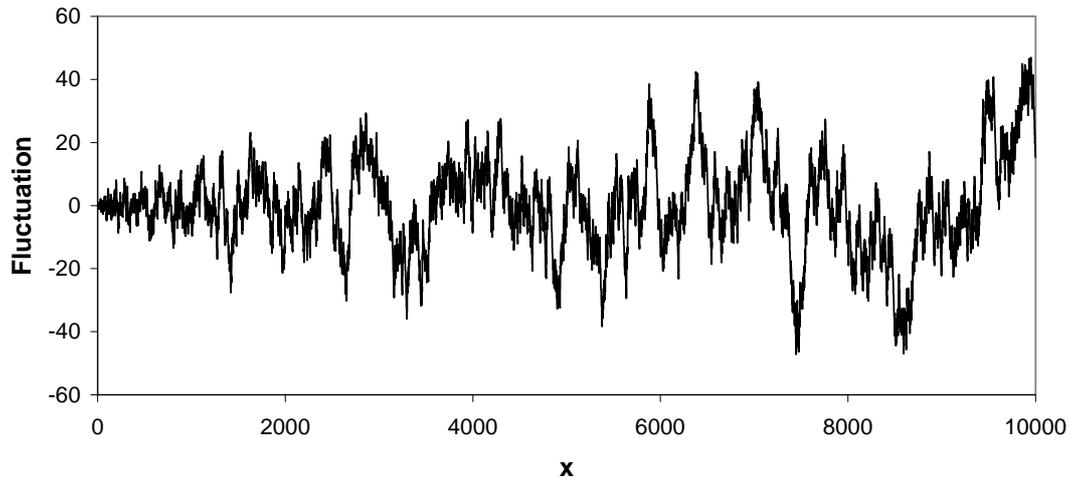

**Fig. 2**

Log-log plot of the one-sided power spectrum of the fluctuating function $\psi_{fluc}(x)$ that was estimated based on $10^3$, $10^4$, $10^5$, and $10^6$ sample points. The power spectrum shifts upwards as the sample size increases.

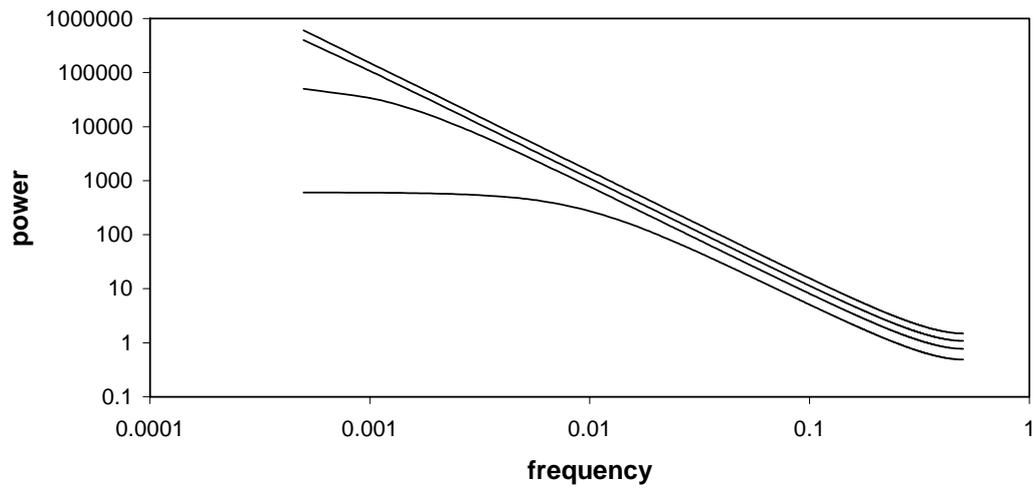

**Fig. 3**

Log-log plot of the one-sided power spectrum (bold line) of the fluctuating function $\psi_{fluc}(x)$, which was estimated based on $10^6$ sample points, together with the straight-line fit (thin line).

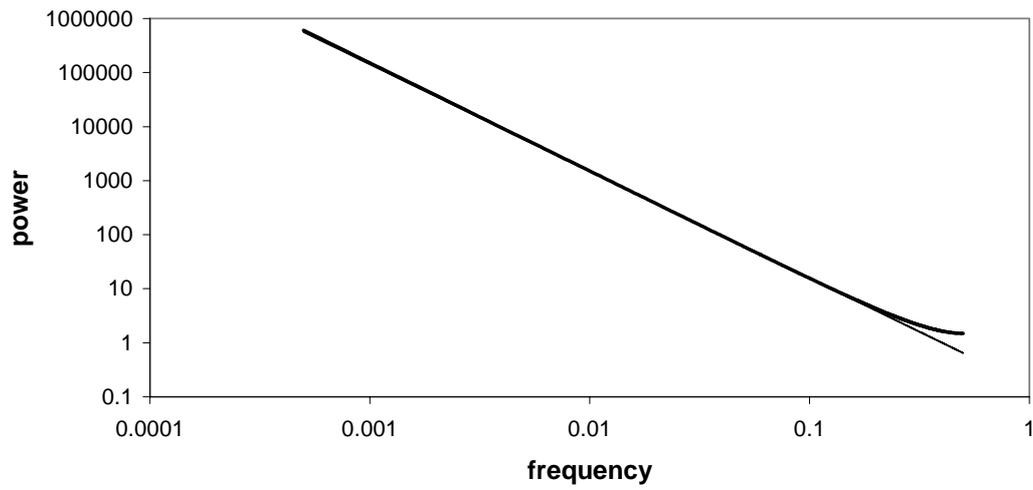